\documentclass{article}
\usepackage[english]{babel}
\usepackage{amssymb,amsmath,amsthm}

\theoremstyle{plain}
\newtheorem{proposition}{Proposition}

\newtheorem{theorem}[proposition]{Theorem}
\newtheorem{lemma}[proposition]{Lemma}

\theoremstyle{definition}
\newtheorem{definition}[proposition]{Definition}

\theoremstyle{remark}
\newtheorem*{remark}{Remark}

\newcommand{\abs}[1]{\left\vert#1\right\vert}
\newcommand{\R}{\mathbb{R}}
\newcommand{\lie}[1]{\mathfrak{#1}}     
\newcommand{\Lie}{\mathcal{L}}          

\newcommand{\hook}{\lrcorner\,}

\newcommand{\SU}{\mathrm{SU}}
\newcommand{\SO}{\mathrm{SO}}

\newcommand{\dfn}[1]{\emph{#1}}

\DeclareMathOperator{\tr}{tr}

\DeclareMathOperator{\ad}{ad}

\newcommand{\Span}[1]{\operatorname{Span}\left\{#1\right\}}

\begin{document}
\title{Solvable Lie algebras are not that hypo}
\author{Diego Conti, Marisa Fern\'andez and Jos\'e A. Santisteban}
\maketitle

\begin{abstract}
We study a type of left-invariant structure on Lie groups, or equivalently on Lie algebras.
We introduce obstructions to the existence of a
hypo structure, namely the $5$-dimensional
geometry of hypersurfaces in manifolds with holonomy $\SU(3)$.
The choice of a splitting $\lie{g}^*=V_1\oplus V_2$, and the vanishing of certain associated
cohomology groups, determine a first obstruction. We also construct necessary conditions for the existence of a hypo structure with a fixed almost-contact form. For non-unimodular Lie algebras, we derive an obstruction to the existence of a hypo structure, with no choice involved.
We apply these methods to classify solvable Lie algebras that admit a hypo structure.
\end{abstract}

\vskip5pt\centerline{\small\textbf{MSC classification}: Primary 53C25;
Secondary 53C15, 17B30, 53D15}\vskip15pt

\section*{Introduction}

In \cite{ContiSalamon}, Salamon and the first author of the present work introduced
{\it hypo structures}, namely the $\SU(2)$-structures induced naturally on orientable hypersurfaces
of Calabi-Yau manifolds of (real) dimension 6.
They are defined as follows.
An $\SU(2)$-structure on a five-manifold is an almost-contact metric structure with additionally a reduction from the structure group $\SO(4)$ to $\SU(2)$; such a structure is entirely determined by the choice of differential forms $(\alpha,\omega_1,\omega_2,\omega_3)$, where $\alpha$ is the almost-contact $1$-form and the $\omega_i$ are pointwise a distinguished orthonormal basis of $\Lambda^2_+(\ker\alpha)$,
which implies that the quadruplet $(\alpha,\omega_1,\omega_2,\omega_3)$ satisfies
certain relations (see Section \ref{sec:coherent}, Proposition \ref{prop:quadruplet}).
Since $\SU(2)$ is the stabilizer of a point under the action of $\SU(3)$ on $\R^6$, hypersurfaces in manifolds
 with holonomy contained in $\SU(3)$  or,
equivalently, with an integrable $\SU(3)$-structure, inherit a natural $\SU(2)$-structure.

In fact, if $M$ is a Riemannian 6-manifold with holonomy contained in
$\SU(3)$, then $M$ has a Hermitian structure, with K\"ahler form $F$, and
a complex volume form $\Psi = \Psi_{+} + i \Psi_{-}$, satisfying $dF=0=d\Psi$. Therefore, if $N\subset M$ is an orientable
hypersurface, and $\mathbb U$ is the unit normal
vector field, the $\SU(3)$-structure induces an
$\SU(2)$-structure $(\alpha,\omega_1,\omega_2,\omega_3)$ on $N$
defined by
  $$
  \alpha=-\mathbb U\lrcorner F,\quad
  \omega_1=f^*F,\quad \omega_2=\mathbb U\lrcorner\Psi_-, \quad
  \omega_3=-\mathbb U\lrcorner\Psi_+,$$ where $f\colon N\to M$ is the inclusion.

The integrability condition on the ambient manifold immediately gives
\[d\omega_1=0, \quad d(\omega_2\wedge\alpha)=0, \quad d(\omega_3\wedge\alpha)=0.\]
An $\SU(2)$-structure satisfying this condition is called a \dfn{hypo} structure. Such a structure can also be characterized in terms of generalized Killing spinors, or  by the condition that the intrinsic torsion is a symmetric tensor, which turns out to coincide with the second fundamental form of the hypersurface. In this sense, hypo geometry is the five-dimensional analogue of half-flat geometry in dimension six (see \cite{ChiossiSalamon, Hitchin:StableForms}). Indeed, much as in the half-flat case, any real-analytic hypo manifold can be immersed isometrically in a Riemannian manifold with holonomy contained in $\SU(3)$, so as to invert the construction outlined above, and the immersion can be determined explicitly by solving certain evolution equations (\cite{ContiSalamon}).

In order to construct examples of hypo structures, a natural place to look is left-invariant structures
on $5$-dimensional
Lie groups. In the analogous half-flat case, this was the approach of \cite{ChiossiSwann,ChiossiFino,ContiTomassini, Conti:Halfflat}, focusing on the nilpotent case, and more recently of \cite{Hengesbach}, considering products of three-dimensional Lie groups. In five-dimensions, only $9$ isomorphism classes of nilpotent Lie groups exist, of which exactly six admit a hypo structure \cite{ContiSalamon}. If one considers solvable Lie groups, things become more complicated. By Mubarakzyanov's classification \cite{Mubarakzjanov}, there are $66$ families of solvable Lie algebras of dimension $5$, some of which depend on parameters; we refer to the comprehensive list of \cite{Basarab-Horwath}. It was shown in \cite{Diatta} that precisely $35$
out of these $66$ families admit an invariant contact structure, at least generically (i.e. for generic values of the parameters). Moreover,  without using Mubarakzyanov's classification,
it was proved in \cite{DeAndres:HypoContact}
that only $5$ of the $66$ admit a hypo-contact structure, namely a hypo structure $(\alpha,\omega_i)$ such that the underlying almost-contact metric structure is a contact metric structure.

\medskip
In this paper we introduce some obstructions to the existence of a hypo structure on a Lie algebra, and use them to classify solvable Lie algebras with a hypo structure. The first obstruction follows a construction of \cite{Conti:Halfflat}. One considers a splitting $\lie{g}^*=V_1\oplus V_2$, where $V_1$ has dimension two. This determines a doubly graded vector space $\Lambda^*\lie{g}^*=\bigoplus\Lambda^{p,q}$, which is made into a double complex if
\begin{equation}
 \label{eqn:coherence}
 d(\Lambda^{p,q})\subset \Lambda^{p+2,q-1}\oplus\Lambda^{p+1,q}.
\end{equation}
The double complex has an associated spectral sequence that collapses at the second step. If $H^{0,3}=E_2^{0,3}$ and $H^{0,2}=E_2^{0,2}$ are zero, relative to some choice of the splitting, then no hypo structure exists (see Proposition $3$). In fact, the key property is
\[Z^k\subset\Lambda^{2,k-2}\oplus\Lambda^{1,k-1},\quad k=2,3,\]
where $Z^k$ denotes the space of closed $k$-forms; this condition does not require \eqref{eqn:coherence}, whose main relevance is in giving a cohomological interpretation. This obstruction applies to $27$
 indecomposable Lie algebras
and $10$
decomposable Lie algebras,
at least generically, where decomposable means isomorphic to a direct sum of ideals.

A second set of obstructions comes from the fact that if $(\alpha,\omega_i)$ is a hypo structure
on a Lie algebra $\lie{g}$, then the forms $\omega_2\wedge\alpha$, $\omega_3\wedge\alpha$ lie in the space
\[V=\{\gamma\in \Lambda^3\lie{g}^*\mid \gamma\wedge\alpha=0, d\gamma=0\}.\]
If, for some $\beta\in\lie{g}^*$, either the space  $V\wedge\beta\subset\Lambda^4\lie{g}^*$ has dimension one or
\begin{equation}
 \label{eqn:secondobstr}
 \dim (V\wedge \beta)=2, \quad Z^2\wedge\alpha\wedge\beta\subset V\wedge\beta,
\end{equation}
then necessarily $\alpha$ and $\beta$ are linearly dependent. This is an obstruction to the existence of a hypo structure with a fixed $\alpha$ (see Proposition $4$), but it can be combined with other arguments to prove that no hypo structure exists on
a Lie algebra.

Indeed, we show that if a non-unimodular Lie algebra $\lie{g}$ has a hypo structure $(\alpha,\omega_i)$, then the $1$-form $\beta\in\lie{g}^*$ defined by $\beta(X)=\tr\ad(X)$ is orthogonal to $\alpha$; this gives a canonical choice for $\beta$ in \eqref{eqn:secondobstr}. Explicitly,
in Proposition $6$,
we prove that there is no hypo structure if either $Z^3\wedge\beta$ has dimension less than two, or
\[\dim (Z^3\wedge \beta)=2, \text{ and }  \alpha\wedge\beta\wedge Z^2\subset Z^3\wedge\beta\]
for any $\alpha\in\lie{g}^*$ such that  $\alpha\wedge\beta\wedge Z^3=0$.
This obstruction applies to $6$ indecomposable families and $12$ decomposable families.

On the other hand, even for unimodular Lie algebras, the structure of the space of closed $3$-forms may give restrictions on $\alpha$ (see Lemma~\ref{lemma:alphaXzero}), which together with \eqref{eqn:secondobstr} enable one to show that certain Lie algebras have no structure. This obstruction accounts for $6$ indecomposable families.

Finally, for $2$ indecomposable families and one decomposable Lie algebra, we use the trivial fact that the space $(Z^2)^2\wedge\alpha$ is non-zero, as it contains $(\omega_1)^2\wedge\alpha\neq 0$.

\medskip
Having obtained the classification, we can ask how often a solvable Lie algebra is hypo. We know from \cite{ContiSalamon} that the answer is $6$ times out of $9$ for nilpotent Lie algebras. In fact, we obtain a shorter proof of this result, namely that the nilpotent Lie algebras denoted here by $D_3$, $A_{5,3}$, $A_{4,1}\oplus\R$ have no hypo structure.

In the solvable case, the question is somewhat ambiguous, because the Lie algebras come in families. With reference to Mubarakzyanov's list, it turns out that, given a family with more than one element, the subset of Lie algebras that have a hypo structure is always a proper subset, but not always discrete. This suggests recasting the question in the following form: how many families in Mubarakzyanov's list of solvable Lie algebras contain at least one hypo Lie algebra? The answer is $21$ out of $66$, so the ratio is considerably less than in the nilpotent case.

If we further distinguish according to whether a family of Lie algebras is decomposable
and whether it is generically contact, we obtain the following table:
\[\begin{array}{l|rr|r|}
  & \text{generically contact} & \text{non-contact} & \text{all}\\
 \hline
 \text{indecomposable} & 7/24 & 9/15 & 16/39\\
 \text{decomposable} & 1/11 & 4/16 & 5/27 \\
 \hline
 \text{all} & 8/35 & 13/31 & 21/66\\
 \hline
\end{array}\]
For instance, the top-left entry states that of the $24$ families in Mubarakzyanov's list which are
indecomposable and have a contact structure for generic choice of the parameters, precisely $7$ have a hypo structure for some choice of the parameters.

If we count Lie algebras with a hypo structure rather than families, we obtain the following table:
\[\begin{array}{l|rr|}
  & \text{generically contact} & \text{non-contact} \\
 \hline
  \text{indecomposable} & 9 &  \text{infinite} \\
 \text{decomposable} & 1 & 4  \\
 \hline
\end{array}\]
Thus, there are exactly ten solvable Lie algebras that have both a hypo and a contact structure, for half of which the structures can be chosen to be compatible  \cite{DeAndres:HypoContact}.

Finally, we point out that there are only five non-unimodular hypo Lie algebras, contained in three families, all of them
indecomposable and contact.

\section{A first obstruction}
\label{sec:coherent}
In this section we introduce an obstruction to the existence of a
hypo structure on a $5$-dimensional Lie algebra. This obstruction
is given in terms of the cohomology groups of a certain double complex
associated to any $n$-dimensional Lie algebra.

Let $\lie{g}$ be an $n$-dimensional Lie algebra, and denote by $d$ the Chevalley-Eilenberg
differential on the dual $\lie{g}^*$. A {\em coherent splitting} of
 $\lie{g}$ is a splitting $\lie{g}^*=V_1\oplus V_2$, where $V_1$ and $V_2$  are vector spaces, $\dim V_1=r \geq 2$ and
\[d(V_1)\subset\Lambda^2V_1, \quad d(V_2)\subset \Lambda^2V_1+V_1\wedge V_2.\]
Let $\Lambda^{p,q}$ be the natural image of $\Lambda^pV_1\otimes \Lambda^qV_2$ in $\Lambda^{p+q}=\Lambda^{p+q}\lie{g}^*$, with the convention that $\Lambda^{p,q}=0$ whenever $p$ or $q$ is negative.
A coherent splitting determines a double
complex $(\Lambda^{*,*}, \delta_1, \delta_2)$,
$\delta_1$, $\delta_2$ being the operators:
\[\delta_1 \colon\Lambda^{p,q}\longrightarrow \Lambda^{p+1,q}, \quad \delta_2 \colon\Lambda^{p,q}\longrightarrow \Lambda^{p+2,q-1}, \quad d=\delta_1+\delta_2.\]
They satisfy
\[\delta_1^2=0=\delta_2^2=\delta_1\delta_2+\delta_2\delta_1.\]
For any choice of coherent splitting on $\lie{g}$, we can define the cohomology
groups $H^{p,q}(\lie{g},V_1)$ as follows (see also \cite{Conti:Halfflat}).
For each $k\geq0$ we define a filtration
\begin{multline}
\label{eqn:Filtration}
\Lambda^{r,k-r}\subset \Lambda^{r,k-r}+\Lambda^{r-1,k-r+1}\subset
\Lambda^{r,k-r}+\Lambda^{r-1,k-r+1}+\Lambda^{r-2,k-r+2}\subset
\\
\cdots
\subset \Lambda^{r,k-r}+\Lambda^{r-1,k-r+1}+\cdots +\Lambda^{0,k}=\Lambda^{k}.
\end{multline}
Notice that in \eqref{eqn:Filtration}, the space $\Lambda^{p,k-p}$ is zero if $p>k$ or $p>r$.
We denote by $Z^k\subset \Lambda^k$ the space of closed invariant $k$-forms. Taking the intersection with $Z^{k}$,  the  filtration \eqref{eqn:Filtration} determines the filtration
\[Z^{k}_r\subset Z^{k}_{r-1}\subset Z^{k}_{r-2}\subset \cdots\subset Z^{k}_0=Z^{k},\]
and taking the quotient by the $d$-exact forms, we obtain yet another filtration
\[H^{k}_r\subset H^{k}_{r-1}\subset H^{k}_{r-2}\subset\cdots\subset H^{k}_0=H^{k}.\]
We can now define the cohomology groups
\[H^{p,q}(\lie{g},V_1)=\frac{H^{p+q}_p}{H^{p+q}_{p+1}}.\]
The notation is justified by the fact that whilst the spaces $\Lambda^{p,q}$ depend on both $V_1$ and $V_2$, the filtration~\eqref{eqn:Filtration}, and therefore the cohomology groups, depend only on
$V_1$. We define
\[h^{p,q}(\lie{g},V_1)=\dim H^{p,q}(\lie{g},V_1).\]

We can think of a coherent splitting as defined by a decomposable form which spans
$\Lambda^{r}V_1$.
\begin{lemma}
\label{lemma:coherentsplitting}
Let $\lie{g}$ be a Lie algebra of dimension $n$, and let $\phi$ be a decomposable $r$-form.
Then $\phi$ defines a coherent splitting $\lie{g}^*=V_1\oplus V_2$, with $\dim V_1=r$,
if and only if
\begin{itemize}
\item $d\alpha\wedge\phi=0$ for all $\alpha\in\lie{g}^*$;
\item $d\phi=0$;
\item $\Lie_X\phi$ is a multiple of $\phi$ for all $X$ in $\lie{g}$, where ${\cal L}$ denotes the Lie derivative.
\end{itemize}
\end{lemma}
\begin{proof}
Given a coherent splitting with $\Lambda^{r,0}$ generated by $\phi$,  we have
\[d\phi\in\Lambda^{r+1,0}=\{0\}\; \text{ and } d\alpha\wedge\phi\in\Lambda^{r+2,0}+
\Lambda^{r+1,1} =\{0\},\; \text{ for } \alpha\in\Lambda^1.\]
Also, since $\phi$ is closed, $\Lie_X\phi=d(X\hook\phi)\in\Lambda^{r,0}$
(where $X\hook\cdot$ denotes the contraction by $X$) which is spanned by $\phi$.

To prove the converse, let $\phi=\alpha^1\wedge\dotsm\wedge\alpha^r$, and complete $\alpha^1,\dotsc,\alpha^r$ to a basis $\alpha^1,\dotsc,\alpha^r,\beta^{1},\dotsc,\beta^{n-r}$. The first condition implies that the image of $d\colon\Lambda^1\to\Lambda^2$ is contained in $\Lambda^{2,0}+\Lambda^{1,1}$. All we need to check in order to have a coherent splitting is that $d\alpha^j$ has type $(2,0)$. Suppose otherwise. Then \[(d\alpha^i)^{1,1}=a^i_{jh}\beta_h\wedge\alpha_j.\]
Now since $\phi$ is closed, $d(X\hook\phi)=\Lie_X\phi$
is a multiple of $\phi$ by hypothesis. So we have that
\begin{multline*}
0=(d(\alpha^1\wedge\dotsb\wedge\hat\alpha^i\wedge\dotsb\wedge\alpha^r))^{r-1,1}=\sum_{j\neq i,h}(-1)^{j+i-1} a^j_{ih}\beta_h\wedge\alpha^1\wedge\dotsb\wedge\hat\alpha^j\wedge\dotsb\wedge\alpha^r\\
+
\sum_{j\neq i,h}a^j_{jh}\beta_h\wedge\alpha^1\wedge\dotsb\wedge\hat\alpha^i\wedge\dotsb\wedge\alpha^r.
\end{multline*}
Hence $a^j_{ih}=0$ for all $i\neq j$ and $h$, and $\sum_{j\neq i}a^j_{jh}=0$ for all $i,h$. Therefore $a^i_{jh}=0$ for all $i,j,h$.
\end{proof}

We introduce the following notation. Let $D_j$ be the annihilator of the kernel of $d\colon\Lambda^j\lie{g^*}\to\Lambda^{j+1}\lie{g^*}$. In other words, if $v^1,\dotsc, v^n$ is a basis of the vector space
$(\Lambda^{j+1})^*$  dual to $\Lambda^{j+1}$,
then $D_j$ is spanned by the $v^j\circ d$. Likewise,
for any $\phi\in\Lambda^k\lie{g}^*$,
let $L_j^\phi$ be the annihilator of the kernel of the map
\[\Lambda^j\lie{g}^*\to\Lambda^{j+k}\lie{g}^*, \quad \alpha\to\alpha\wedge\phi.\]

We can then give a specialized version of the lemma that accounts for the vanishing of certain cohomology groups. In the five-dimensional case we get:
\begin{proposition}
\label{prop:criterion}
Let $\lie{g}$ be a $5$-dimensional Lie algebra. Then $\lie{g}$ has a coherent splitting,
with $\dim V_1=2$ and
$H^{0,2}=0=H^{0,3}$,
if and only if
there exists a nonzero $2$-form $\phi$ such that
\begin{itemize}
\item $\phi\wedge\phi=0$;
\item $d\phi=0$;
\item $\Lie_X\phi$ is a multiple of $\phi$ for all $X$ in $\lie{g}$;
\item $L_2^\phi\subset D_2$, $L_3^\phi\subset D_3$.
\end{itemize}
\end{proposition}
\begin{proof}
Given a coherent splitting, it is clear that exact $k$-forms have no component
in $\Lambda^{0,k}$. Moreover,
the condition $H^{0,k}=0$ is equivalent to $Z^k$ being contained in $\Lambda^{k,0}+\dotsb + \Lambda^{1,k-1}$. Thus, $L_k^\phi\subset D_k$ if and only if $H^{0,k}=0$.

Conversely, a $2$-form $\phi$ such that $\phi\wedge\phi=0$ is decomposable, and therefore determines a splitting. If $\phi$ is as in the hypothesis, the splitting is coherent  because  $L_2^\phi\subset D_2$ implies that closed $2$-forms, and in particular exact $2$-forms, have no component in $\Lambda^{0,2}$; therefore,  $d\alpha\wedge\phi=0$ for all $\alpha\in\lie{g}^*$ and
Lemma~\ref{lemma:coherentsplitting} applies.
\end{proof}

\begin{remark}
In the proof of Proposition~\ref{prop:criterion}, we can suppose
that $\lie{g}$ has a coherent splitting, with $\dim V_1=r \geq 2$, and conclude that
$H^{0,k}=0$ is equivalent to $L_k^\phi\subset D_k$, since this works for any dimension
$n$ of $\lie{g}$ and for all values of $r$. However, we need $r=2$ to have that
a $2$-form $\phi$ is decomposable if and only if $\phi\wedge\phi=0$.
\end{remark}

From now on, given a $5$-dimensional Lie algebra $\lie g$ whose dual is spanned by $\{ e^1,\ldots ,e^5\}$, we will write $e^{ij}= e^i\wedge e^j$,
$e^{ijk}= e^i\wedge e^j\wedge e^k$, and so forth.

The relevance of the above proposition comes from {\em hypo geometry}.
First we recall some facts about $\SU(2)$-structures on a
5-manifold. (For more details, we refer to \cite{ContiSalamon}). Let $N$
be a 5-manifold and let $L(N)$ be the principal bundle of linear
frames on $N$. An {\em $\SU(2)$-structure} on $N$ is an
$\SU(2)$-reduction of $L(N)$. We have the following (see
\cite[Proposition 1]{ContiSalamon}):
\begin{proposition}
\label{prop:quadruplet}
$\SU(2)$-structures on a $5$-manifold $N$ are in
one-to-one correspondence with quadruplets $(\alpha,\omega_1,\omega_2,\omega_3)$,
where $\alpha$ is a $1$-form and $\omega_i$ are $2$-forms on $N$
satisfying
at each point
$$
\omega_i\wedge\omega_j=\delta_{ij}v, \quad
v\wedge\alpha\not=0,
$$
for some $4$-form $v$, and
$$
i_X\omega_1=i_Y\omega_2\Rightarrow \omega_3(X,Y)\ge 0,
$$
where $i_X$ denotes the contraction by $X$.
\end{proposition}

Moreover, we need recall the following definition.
\begin{definition}
An $\SU(2)$-structure $(\alpha,\omega_1,\omega_2,\omega_3)$ on a
$5$-manifold $N$ is said to be {\em hypo} if
  \begin{equation}
d(\omega_2\wedge\alpha)=d(\omega_3\wedge\alpha)=d\omega_1=0.
 \end{equation}
\end{definition}

Therefore, to a choice of a coframe
$f^1,\dotsc, f^5$
on a Lie algebra $\lie{g}$, we associate an $\SU(2)$-structure given by
\begin{equation}
 \label{eqn:adaptedframe}
\alpha=f^5, \quad \omega_1=f^{12}+f^{34}, \quad  \omega_2=f^{13}+f^{42}, \quad  \omega_3=f^{14}+f^{23},
\end{equation}
and it is called {\em hypo} if $\omega_1$, $\omega_2\wedge\alpha$, $\omega_3\wedge\alpha$ are closed.
\begin{definition}
Let $f^1,\dotsc, f^5$ be a coframe on a Lie algebra $\lie{g}$ such that
the quadruplet $(\alpha,\omega_1, \omega_2,\omega_3)$ given by
\eqref{eqn:adaptedframe} defines a hypo structure on $\lie{g}$. Then,
the coframe $f^1,\dotsc, f^5$
is said to be {\em adapted}  to the hypo structure $(\alpha,\omega_1, \omega_2,\omega_3)$.
\end{definition}

\begin{proposition}
\label{prop:h02h03}
If $\lie{g}$ has dimension $5$, and there exists a coherent splitting
$\lie{g}^*=V_1\oplus V_2$ with $\dim V_1=2$ and
$h^{0,3}(\lie{g},V_1)=0=h^{0,2}(\lie{g},V_1)$, then there is no hypo structure.
\end{proposition}
\begin{proof}
Let $(\alpha,\omega_i)$ be a hypo structure,
and let {$\phi$} be a generator of $\Lambda^{2,0}$.
We know that the forms $\omega_1$, $\alpha\wedge\omega_2$, $\alpha\wedge\omega_3$ are closed.
Moreover, because $h^{0,3}(\lie{g},V_1)=0=h^{0,2}(\lie{g},V_1)$, we have
\[\phi\wedge\omega_1=0, \quad \phi\wedge(\alpha\wedge\omega_i)=0.\]
If we decompose the space of two-forms on $\R^5$ according to
\[\Lambda^2\R^5=\alpha\wedge\Lambda^1\R^4\oplus\Lambda^2_+\R^4\oplus\Lambda^2_-\R^4,\]
we find that $\phi$ must lie in $\Lambda^2_-\R^4$. Since $\phi^2=0$, this implies $\phi=0$.
\end{proof}

\begin{remark}
Strictly speaking Proposition~\ref{prop:h02h03} does not use the fact that the splitting is coherent, but only the conditions $Z^3\wedge\phi=0$, $Z^2\wedge\phi=0$; or, in the language of Proposition~\ref{prop:criterion}, the inclusions {$L_2^\phi\subset D_2$, $L_3^\phi\subset D_3$}. Indeed it is sometimes the case that a splitting with this property exists withouth being coherent. Consider for instance the Lie algebra
\[\lie{g}=(e^{13},e^{34},-e^{24},0,0)\]
(by this notation we mean that the dual $\lie{g}^*$ has a fixed basis $e^1,\dotsc, e^5$, such that
\[de^1=e^{13},\;de^2=e^{34},\;de^3=-e^{24},\;de^4=0=de^5). \]
Then $e^{24}$ defines a splitting with $Z^3\wedge\phi=Z^2\wedge\phi=0$, and yet this is not coherent. On the other hand, a different obstruction applies to this case (see Proposition~\ref{prop:unimodular} below).
\end{remark}

\section{Other obstructions}
\label{sec:unimodular}
When looking at  $5$-dimensional solvable Lie algebras, the coherent splitting obstruction,
shown in Proposition  $3$,
is sometimes not sufficient to determine whether a hypo structure exists. In this section we describe two different obstructions that can be used in these cases.

For every $1$-form $\gamma$, let $L_\gamma\colon
\Lambda^j\to \Lambda^{j+1}$ be the map given by
$L_\gamma(\eta)=\gamma\wedge\eta$.
\begin{proposition}
\label{prop:third}
Let $\alpha,\beta$ be linearly independent one-forms on a Lie algebra $\lie{g}$, and set $V=\ker L_\alpha\cap Z^3$. Suppose that either
\begin{itemize}
 \item $\dim L_\beta(V)<2$; or
 \item $\dim L_\beta(V)=2$ and \[L_\alpha(L_\beta(Z^2))\subset L_\beta(V).\]
\end{itemize}
Then there is no hypo structure on $\lie{g}$ of the form $(\alpha,\omega_i)$ (in the sense that its almost-contact form is $\alpha$ itself).
\end{proposition}
\begin{proof}
Suppose for a contradiction that a hypo structure $(\alpha,\omega_i)$ exists, and let $f^1,\dotsc,f^5$ be an adapted coframe. Up to rescaling the metric and up to $\SU(2)$ action, we can assume that
$\beta=f^1+af^5$, with $a$ a constant. Then $\omega_2\wedge\alpha$, $\omega_3\wedge\alpha$ lie in $V$ and
\begin{align*}
L_\beta(V)&\ni L_\beta(\omega_2\wedge\alpha)=f^{1425},&
L_\beta(V)&\ni L_\beta(\omega_3\wedge\alpha)=f^{1235}.
\end{align*}
So $\dim L_\beta(V)\geq 2$, and if equality holds then $L_\beta(V)$ is spanned by $f^{1425}$, $f^{1235}$. But then
\[L_\alpha(L_\beta(Z^2))\ni \alpha\wedge\beta\wedge\omega_1 = f^{5134},\]
which cannot lie in $L_\beta(V)$.
\end{proof}

For non-unimodular Lie algebras, it turns out that we have a canonical choice for $\beta$:
\begin{lemma}
 \label{lemma:unimodular}
Let $\lie{g}$ be a non-unimodular Lie algebra and let $\beta\in\lie{g}^*$ be the form corresponding to the linear map $\lie{g}\to\R$, $X\to \operatorname{tr}\ad X$. If $\lie{g}$ has a hypo structure $(\alpha,\omega_i)$, then $\alpha$ and $\beta$ are orthogonal with respect to the underlying metric.
\end{lemma}
\begin{proof}
In an adapted coframe $f^1,\dotsc, f^5$, with dual frame $f_1,\dotsc, f_5$,
$\alpha=f^5$, we have
\[df^{1234}=\sum_k f^k([f_k,f_5]) f^{12345}=-\beta(f_5) f^{12345}.\]
However, since $\omega_1$ is closed, the left-hand side is zero and so $\beta(f_5)=0$.
\end{proof}

Thus, in the non-unimodular case Proposition~\ref{prop:third} gives a fairly straightforward criterion:
\begin{proposition}
\label{prop:unimodular}
Let $\lie{g}$ be a non-unimodular Lie algebra, and let $\beta\in\lie{g}^*$ be the form corresponding to the linear map $\lie{g}\to\R$, $X\to \operatorname{tr}\ad X$. Suppose that either
\begin{itemize}
 \item $\dim L_\beta(Z^3)<2$; or
 \item $\dim L_\beta(Z^3)=2$ and for every $\alpha\in\lie{g}^*$ such that $L_\alpha(L_\beta(Z^3))=0$, \[L_\alpha(L_\beta(Z^2))\subset L_\beta(Z^3).\]
\end{itemize}
Then $\lie{g}$ has no hypo structure.
\end{proposition}
\begin{proof}
Suppose $\lie{g}$ has a hypo structure $(\alpha,\omega_i)$.
By Lemma \ref{lemma:unimodular}, we know that $\alpha$ and $\beta$ are linearly
independent. Consider the space $V=\ker L_\alpha\cap Z^3$.
If $L_\beta(Z^3)$ has dimension two, then $L_\beta(V)\subseteq L_\beta(Z^3)$ may only have dimension two if equality holds, implying $L_\alpha(L_\beta(Z^3))=0$. Then the statement follows from Proposition~\ref{prop:third}.
\end{proof}

In order to apply Proposition~\ref{prop:third} effectively, one needs information on what the $1$-form $\alpha$ can be. The condition of Lemma~\ref{lemma:unimodular} is often  useful but not always sufficient, since in practice it only tells us that $\alpha$ and $\beta$ are linearly independent; moreover, it does not apply to unimodular Lie algebras (for which $\beta$ is zero). The following result gives useful restrictions on the $1$-form $\alpha$; it is labeled a lemma because we view it as an aid toward the application of either Proposition~\ref{prop:third} or Proposition~\ref{prop:unimodular}.
\begin{lemma}
\label{lemma:alphaXzero}
Suppose $\lie{g}$ has a hypo structure $(\alpha,\omega_i)$. If $\dim (X\hook Z^3)\wedge \gamma<2$, where $X\in\lie{g}$ and $\gamma\in\lie{g}^*$, then $\alpha(X)=0$.
\end{lemma}
\begin{proof}
Suppose otherwise; fix an adapted coframe $f^1,\dotsc, f^5$ with dual frame $f_1,\dotsc, f_5$. Then $X=a f_1 + f_5$ up to a multiple and $\SU(2)$ action. Therefore $X\hook Z^3$ contains
\begin{align*}
X\hook (f^{135}+f^{425})&=f^{13}+f^{42} + af^{35},\\
X\hook (f^{145}+f^{235})&=f^{14}+f^{23} + af^{45}.
\end{align*}
Now by hypothesis some linear combination \[\delta=\lambda (f^{13}+f^{42} + af^{35}) + \mu(f^{14}+f^{23} + af^{45})\] gives zero on wedging with $\gamma$. But
\[(\lambda (f^{13}+f^{42} + af^{35}) + \mu(f^{14}+f^{23} + af^{45}))^2 = 2(\lambda^2+\mu^2)f^{1234} + \text{other terms},\]
which is nonzero. By  non-degeneracy $\delta\wedge\gamma\neq 0$, which is absurd.
\end{proof}
\begin{remark}
Regardless of Proposition~\ref{prop:third}, this lemma may have more immediate applications. Indeed, a hypo structure $(\alpha,\omega_i)$ always satisfies
\begin{equation}
 \label{eqn:Z2wedgealpha}
 0\neq (\omega_1)^2\wedge\alpha\in (Z^2)^2\wedge\alpha.
\end{equation}
\end{remark}

\section{Diatta's algebras}
Let us recall firstly that a {\em contact form} $\eta$ on a five-dimensional Lie algebra $\lie{g}$
is a 1-form on $\lie{g}$ (that is, $\eta\in \lie{g}^*$) such that
$$\eta\wedge\ (d\eta)^2 \not=0.$$
The existence of a hypo structure on $\lie{g}$ is independent of the existence of a
contact form. In fact, in this section we will consider indecomposable solvable Lie algebras
of dimension 5 having a contact form $\eta$, and we will see that many of those Lie
algebras do not admit a hypo structure. Notice that we are not requiring that the almost-contact 1-form $\alpha$ associated to the hypo structure coincide with the contact form $\eta$.

In \cite{Diatta}, Diatta gives a list of 24 (families of)
indecomposable five-dimensional solvable Lie algebras $D_1,\dotsc ,D_{24}$ that admit a left-invariant contact $1$-form. They correspond to the algebras $A_{5,k}$ of \cite{Basarab-Horwath} under
\[D_k\to \begin{cases} A_{5,k+3} & k=1,2,3\\ A_{5,k+15}
&4\leq k \leq 24.\end{cases}.\]
We shall use the notation $D_k(p_1,\dotsc, p_n)$ to denote special instances of a family for assigned values of the parameters. Notice that Diatta's list, as well as the one in \cite{Basarab-Horwath} from which it was extracted, contains conditions on the parameters. We shall ignore these conditions to keep things simpler. This has two consequences: first, the same Lie algebra may appear more than once, and second, a Lie algebra $D_k(p_1,\dotsc, p_n)$ may not have a contact structure for some choice of the parameters $p_1,\dotsc, p_n$. However, these ``degenerate'' cases turn out to never have a hypo structure.

\begin{proposition}
\label{prop:diattahypo}
The indecomposable
solvable Lie algebras that have both a contact structure and a hypo structure are the following:
\begin{align*}
  D_1&=(e^{24}+e^{35},0,0,0,0)\\
  D_2&=\left(e^{34}+e^{25},e^{35},0,0,0\right)\\
  D_4(-1/2,-3/2)&=\left(-\frac{1}{2} e^{15}-e^{23},- e^{25},\frac{1}{2}  e^{35},\frac{3}{2}  e^{45},0\right)\\
  D_4(1,-3)&=\left(-e^{23}-2 e^{15},- e^{25},- e^{35},3  e^{45},0\right)\\
  D_4(-2,3)&=\left(e^{15}-e^{23},- e^{25},2  e^{35},-3  e^{45},0\right)\\
  D_{15}(-1)&=\left(-e^{15}-e^{24},- e^{34},e^{35},- e^{45},0\right)\\
  D_{18}(-1,-1)&=(-e^{14}, -e^{25}, e^{34}+e^{35},0,0)\\
  D_{20}(-2,0)&=(2e^{14},-e^{24}-e^{35},e^{25}-e^{34},0,0)\\
  D_{22}&=\left(e^{23}+2e^{14},e^{24}+e^{35},e^{34}-e^{25},0,0\right)
\end{align*}
A Lie algebra $D_k(p_1,\dotsc, p_n)$ is hypo if and only if it belongs to this list.
\end{proposition}
\begin{proof}
First, we produce a hypo structure for each Lie algebra appearing in the statement.

The Lie algebras $D_4(1,-3)$, $D_{15}(-1)$,  $D_{18}(-1,-1)$ and $D_{22}$ appear in \cite{DeAndres:HypoContact} and have hypo contact structures given by the coframes
\begin{align*}
&e^5,\frac15(e^1-e^4),\frac12e^3,\frac15e^2,-\frac15e^1-\frac2{15}e^4  & D_4&(1,-3). \\
&\frac12(-e^1+e^3),e^5,  \frac{\sqrt2}2e^4,  \frac{\sqrt2}2e^2, -e^1-e^3 & D_{15}&(-1).\\
&-\frac1{2\sqrt3}(e^4+2e^5),\frac1{2\sqrt3}(e^2-e^3),\frac13e^1-\frac16e^2-\frac16e^3,
\frac12e^4,\frac13(e^1+e^2+e^3)& D_{18}&(-1,-1).\\
&e^4,e^1,\frac{\sqrt2}2e^3,\frac{\sqrt2}2e^2,\frac13(e^5-3e^1) & D_{22}&.
\intertext{The Lie algebra $D_1$  is nilpotent and has a well-known hypo-contact structure:}
&e^2,e^4,e^3,e^5,e^1 & D_1&.
\intertext{The Lie algebra $D_2$ is also nilpotent and equivalent to $(0,0,0,12,13+24)$, hence hypo by \cite{ContiSalamon}; a hypo structure is given by the
coframe}
&e^{2},e^{4},e^{5},- e^{1},-e^{1}+e^{3} & D_2&.
\intertext{Hypo structures on the three remaining Lie algebras are new. They are defined by}
& e^{1},e^{3},e^{2},e^{5},e^{4} & D_4&(-1/2,-3/2).\\
&- e^{3},2  e^{5},-2  e^{1},2  e^{2},- e^{4} \sqrt{2} & D_4&(-2,3).
\end{align*}
and, for $D_{20}(-2,0)$, by
\begin{multline*}
 3  e^{2},-3  \sqrt{3} e^{1}- \sqrt{3} e^{4}+2  \sqrt{3} e^{5}-2  \sqrt{3} e^{2}- \sqrt{3} e^{3},9 e^{1}+3 e^{3}+3 e^{4},\\
 -2  \sqrt{3} e^{4}- \sqrt{3} e^{2},-5 e^{1}-2 e^{2}-4 e^{3}-e^{4}+2 e^{5}.
 \end{multline*}
It is straightforward to check that all these Lie algebras have a contact form.

It remains to show that the remaining $D_k(p_1,\dotsc, p_n)$ do not have a hypo structure; to that effect, we apply the results of Sections~\ref{sec:coherent},~\ref{sec:unimodular}. Looking at the list and applying Proposition~\ref{prop:criterion}, we see that the algebras in the list that admit a coherent splitting with \[H^{0,2}=H^{0,3}=0\] are precisely the following ($\phi$ denotes a generator of $\Lambda^{2,0}$ in each case):
\begin{itemize}
\item $D_4=\left(- {(1+p)} e^{15}-e^{23},- e^{25},- p e^{35},- q e^{45},0\right)$, $\phi=e^{25}$
if all of  $p+q$, $2p+1$, $2p+1+q$ are non-zero, or $\phi=e^{35}$ if all of $p+q+2$, $p+2$, $1+q$ are non-zero;
\item $D_5=\left(- e^{15} {(1+p)}-e^{23}-e^{45},- e^{25},- p e^{35},- e^{45} {(1+p)},0\right)$, $\phi= e^{25}$
if both of $1+2p$, $2+3p$ are non-zero, or $\phi= e^{35}$ otherwise;
\item $D_6=\left(2 e^{15}+e^{23},e^{25},e^{25}+e^{35},e^{45}+e^{35},0\right)$, $\phi=e^{25}$;
\item $D_7=\left(e^{23},0,e^{25},e^{45},0\right)$, $\phi=e^{25}$;
\item $D_8=\left(-2 e^{15}-e^{23},- e^{25},-e^{25}-e^{35},- p e^{45},0\right)$, $\phi=e^{25}$;
\item $D_9=\left(2 e^{15}+e^{23}+\epsilon e^{45},e^{25},e^{35}+e^{25},2  e^{45},0\right)$, $\phi=e^{25}$;
\item $D_{12}=\left(e^{45}+e^{15}+e^{23},0,e^{35},e^{35}+e^{45},0\right)$, $\phi= e^{25}$;
\item $D_{13}=\left(- e^{15} {(1+p)}-e^{23},- p e^{25},-e^{35}- p e^{25},-e^{35}-e^{45},0\right)$, $\phi=e^{25}$
if both of $p+2$, $p+3$ are non-zero, or $\phi=p e^{25}+(1-p)e^{35}$ otherwise;
\item $D_{14}=\left(e^{15}+e^{23},e^{25},0,e^{45},0\right)$, $\phi=e^{25}$;
\item $D_{15}=\left(-e^{24}- e^{15} {(2+p)},- e^{25} {(1+p)}-e^{34},- p e^{35},- e^{45},0\right)$, $\phi=e^{45}$
if both of $1+2p$, $1+p$ are non-zero;
\item $D_{16}=\left(e^{24}+3 e^{15},e^{34}+2 e^{25},e^{35}+e^{45},e^{45},0\right)$, $\phi=e^{45}$;
\item $D_{17}=\left(-e^{15}-e^{24}- p e^{35},-e^{25}-e^{34},- e^{35},0,0\right)$, $\phi=e^{45}$;
\item $D_{18}=\left(- e^{14},- e^{25},- p e^{34}- q e^{35},0,0\right)$, $\phi=e^{45}$ if $(p,q)$ is not $(0,-1)$, $(-1,0)$ or $(-1,-1)$;
\item $D_{19}=\left(-e^{15}- p e^{14},-e^{35}-e^{24},- e^{34},0,0\right)$, $\phi=e^{45}$'
\item $D_{20}=\left(- q e^{15}- p e^{14},-e^{24}-e^{35},e^{25}-e^{34},0,0\right)$, $\phi=e^{45}$ if $(p,q)$ is not $(-2,0)$;
\item $D_{23}=\left(e^{14},e^{25},e^{45},0,0\right)$, $\phi=e^{45}$;
\item $D_{24}=\left(e^{14}+e^{25},e^{24}-e^{15},e^{45},0,0\right)$, $\phi=e^{45}$.
\end{itemize}

Other non-unimodular Lie algebras are ruled out by Proposition~\ref{prop:unimodular}. They are listed below; here and throughout the paper, the $1$-form $\beta$ is given up to multiple.
\begin{itemize}
\item $D_4(1, -1)=\left(-e^{23}-2 e^{15},- e^{25},- e^{35},e^{45},0\right)$, $\beta=e^5$
\item $D_4(0, -1)=\left(-e^{15}-e^{23},- e^{25},0,e^{45},0\right)$, $\beta=e^5$
\item $D_{15}(-1/2)=\left(-e^{24}-\frac{3}{2} e^{15},-\frac{1}{2} e^{25}-e^{34},\frac{1}{2}  e^{35},- e^{45},0\right)$, $\beta=e^5$.
\end{itemize}

To address the remaining Lie algebras, we apply either Proposition~\ref{prop:third} or Equation~\ref{eqn:Z2wedgealpha}.
\begin{itemize}
\item $D_3=\left(e^{25}+e^{34},e^{35},e^{45},0,0\right)$. This Lie algebra is nilpotent and isomorphic to $(0,0,12,13,23+14)$, therefore not hypo by \cite{ContiSalamon}; however, we can prove it directly using the methods of Section~\ref{sec:unimodular}. We compute
\[Z^2=\operatorname{Span}\{-e^{14}+e^{23},e^{15}+e^{24},e^{25},e^{34},e^{35},e^{45}\}\]
\[Z^3=\operatorname{Span}\{-e^{125}+e^{134},e^{135},e^{145},e^{234},e^{235},e^{245},e^{345}\}\]
Then the spaces $(e_i\hook Z^3)\wedge e^5$, $i=1,2$ are one-dimensional, so by Lemma~\ref{lemma:alphaXzero} $\alpha$ is a linear combination of $e^3,e^4,e^5$. In particular, if $\alpha$ is linearly independent of $e^5$, then setting $\beta=e^5$ in Proposition~\ref{prop:third} we see that
\[L_\beta(V)\subset\Span {e^{1345},e^{2345}}\] contains $Z^2\wedge\alpha\wedge\beta$. Otherwise, we may set $\beta=e^4$ and obtain the same result.
\item $D_4(-2,2)=\left(e^{15}-e^{23},- e^{25},2  e^{35},-2  e^{45},0\right)$. We compute
\[
Z^2=\Span{e^{12},-e^{15}+e^{23},e^{25},e^{34},e^{35},e^{45}}\]
\[
Z^3=\Span{e^{125},
e^{135},
e^{235},e^{234}-e^{145},
e^{245},
e^{345}}\]
Therefore $Z^3\wedge e^5$ is one-dimensional, hence by Lemma~\ref{lemma:alphaXzero} $\alpha(e_i)=0$, $i=1,2,3,4$, i.e $\alpha=e^5$ up to a multiple. Now $Z^2\wedge\alpha$ is spanned by $e^{125},e^{235},e^{345}$. Setting $\beta=e^3$ in Proposition~\ref{prop:third} gives a contradiction, as $L_\beta(V)$ is spanned by $e^{1235},e^{2345}$ and $Z^2\wedge\alpha\wedge\beta$ is spanned by $e^{1235}$.
\item $D_4(-1/2,-1)=\left(-\frac{1}{2} e^{15}-e^{23},- e^{25},\frac{1}{2}  e^{35},e^{45},0\right)$ is similar to $D_4(-2,2)$ in that
\[Z^2=\Span{e^{13},
\frac{1}{2} e^{15}+e^{23},e^{24},e^{25},e^{35},e^{45}},\]
\[Z^3=\Span{e^{125},
e^{135},e^{234}+\frac{1}{2} e^{145},
e^{235},
e^{245},
e^{345}}.\]
The same argument applies, except that now $Z^2\wedge\alpha\wedge\beta$ is spanned by $e^{2345}$.
\item $D_{10}=\left(-2  p e^{15}-e^{23},- p e^{25}+e^{35},- p e^{35}-e^{25},- q e^{45},0\right)$.
A basis of $Z^3$ is given by
\[
 e^{345}, e^{235} ,e^{135}, e^{125}, e^{245}, (2p+q)e^{145} - e^{234},
\]
plus  $e^{123}$  if $p=0$, whereas
\[Z^2=\Span{e^{35}, e^{45}, 2 e^{15}+  e^{23}, e^{25}}.\]
Now $\beta=(4p+q)e^5$, and if $4p+q\neq 0$ then Proposition~\ref{prop:unimodular} applies; in general, we have that $e_1\hook Z^3$ and $e_4\hook Z^3$ wedged with $e^5$ are at most one-dimensional, hence by Lemma~\ref{lemma:alphaXzero} $\alpha$ lies in the span of $e^2,e^3,e^5$. But then $(Z^2)^2\wedge\alpha=0$, which is a contradiction.
\item $D_{11}=\left(-2  e^{15} p-e^{23}-\epsilon e^{45},- e^{25} p+e^{35},- e^{35} p-e^{25},-2  e^{45} p,0\right)$.
A basis of $Z^2$ is given by
\[ 2p e^{15}+ e^{23}, e^{35}, e^{45}, e^{25},\]
whereas a basis of $Z^3$ is given by
\[ e^{345},e^{235}, e^{135},e^{125}, e^{245},-4p e^{145}+ e^{234},\]
plus $e^{145}-\epsilon e^{123} $ if $p=0$.
Thus $(e_1\hook Z^3)\wedge e^5$, $(e_4\hook Z^3)\wedge e^5$ have dimension one and we see that $\alpha$ is in $\Span{e^2,e^3,e^5}$, whence $\alpha\wedge (Z^2)^2=0$, a contradiction.
\item $D_{18}(-1,0)=(-e^{14}, -e^{25}, e^{34}, 0, 0)$. We compute
\[Z^2=\Span{e^{13},e^{14},e^{25},e^{34},e^{45}},\]
\[Z^3=\Span{e^{125}+e^{124},e^{134},e^{135},e^{145},-e^{234}+e^{235},e^{245},e^{345}},\]
and $\beta=e^5$. Moreover the spaces \[(e_1\hook Z^3)\wedge (e^4+e^5),\quad (e_2\hook Z^3)\wedge (e^4+e^5),\quad (e_3\hook Z^3)\wedge (e^4-e^5)\] are one-dimensional, so by Lemma~\ref{lemma:alphaXzero} and Lemma~\ref{lemma:unimodular} we get $\alpha=e^4+ae^5$, for some constant $a$. Then setting $\beta=e^5$ in Proposition~\ref{prop:third}, we see that $L_\beta(V)$ is at most two-dimensional, and it contains $e^{1345}$. Since $Z^2\wedge\alpha\wedge\beta$ is spanned by $e^{1345}$, there is no hypo structure.
\item $D_{18}(0,-1)=(-e^{14}, -e^{25}, e^{35}, 0, 0)$ is really isomorphic to $D_{18}(-1,0)$,
as one can check by considering the coframe $(e^2,e^1,e^3,e^5,e^4)$, so it has no hypo structure.
\item $D_{21}=(e^{23}+e^{14},e^{24}-e^{25},e^{35},0,0)$. Then $\beta=e^4$ and
\begin{align*}
 Z^2&=\Span{e^{14}+e^{23},e^{25}-e^{24},e^{35},e^{45}},\\
Z^3&=\Span{e^{125}-2 e^{124},e^{135}+e^{134},e^{234},e^{235}+e^{145},e^{245},e^{345}}.
\end{align*}
Therefore, the spaces \[(e_1\hook Z^3)\wedge (e^4+e^5), \quad (e_2\hook Z^3)\wedge (e^5-2e^4), \quad (e_3\hook Z^3)\wedge (e^4+e^5)\] are one-dimensional, so by Lemma~\ref{lemma:alphaXzero} and Lemma~\ref{lemma:unimodular} we get $\alpha=ae^4+e^5$, for some constant $a$.  Then setting $\beta=e^4$ in Proposition~\ref{prop:third}, we see that $L_\beta(V)$ is at most two-dimensional, and it contains $e^{2345}$. Since $Z^2\wedge\alpha\wedge\beta$ is spanned by $e^{2345}$, there is no hypo structure.
\end{itemize}
\end{proof}

\section{Indecomposable Lie algebras without contact form}
We now pass on to indecomposable solvable Lie algebras that do not have a contact structure.
\begin{proposition}
The indecomposable
solvable Lie algebras which have a hypo structure but not a contact structure are those given in Table~\ref{table:IrreducibleNonContact}, all of them unimodular.
\end{proposition}
Observe that $A_{5,1}$ and $A_{5,2}$ are nilpotent, and so appear in \cite{ContiSalamon}.
\begin{proof}
It is straightforward to verify that the coframes given in the table define indeed hypo structures. To show that no other Lie algebras of the specified type have a hypo structure, we use the classification in \cite{Basarab-Horwath}.
\begin{itemize}
\item $A_{5,3}=\left(e^{25},e^{45},e^{24},0,0\right)$ is nilpotent and
known not to have a hypo structure \cite{ContiSalamon}. It also has a coherent splitting $\phi=e^{45}$ with $H^{0,2}=H^{0,3}=0$.

\item $A_{5,7}=\left(e^{15}, p e^{25}, qe^{35} , r e^{45},0\right)$ where $p,q,r\neq0$ is not hypo unless, up to permutation of the parameters, $r=-1$ and $p+q=0$. Indeed, suppose first that $p,q,r\neq-1$. Then  if $p+q\neq0,-1$ we find a coherent splitting $\phi=e^{45}$ with $H^{0,2}=H^{0,3}=0$. Since we can act by an automorphism to permute $p,q,r$, the same happens if $p+r\neq0,-1$ or $q+r\neq0,-1$. Thus, still assuming $p,q,r\neq-1$, we are left with the cases $(-\frac12,-\frac12,-\frac12)$ and $(-\frac12,-\frac12,\frac12)$. In the former case, $\phi=e^{15}$ defines a coherent splitting with $H^{0,2}=H^{0,3}=0$. In the latter case, the Lie algebra is non-unimodular with $\beta=e^5$, and $L_\beta(Z^3)$ has dimension one.

Now, if $r=-1$ and $p+q\neq0$, the Lie algebra is non-unimodular with $\beta=e^5$, and $Z^3$ is spanned by
\[  e^{345}, e^{145}, e^{235}, e^{125}, e^{245}, e^{135}\]
plus
\begin{gather*}
 e^{234} \text{ if } p+q-1=0,\\
 e^{123} \text{ if } p+q+1=0.
\end{gather*}
Therefore,  $L_\beta(Z^3)$ is at most one-dimensional.
\item $A_{5,8}=\left(e^{25},0,e^{35}, pe^{45},0\right)$ has a coherent splitting  given by $\phi=e^{25}$ with $H^{0,2}=H^{0,3}=0$ if $p\neq -1$.
\item $A_{5,9}=\left(e^{15}+e^{25},e^{25},pe^{35}, qe^{45} ,0\right)$ has a coherent splitting with
$H^{0,2}=H^{0,3}=0$. If $p+q\neq 0$ we can take $\phi=e^{25}$; if $p=-q$ but $p\neq 1,2$ then $\phi=e^{35}$, and otherwise we can take $\phi=e^{45}$.
\item $A_{5,10}=\left(e^{25},e^{35},0,e^{45},0\right)$ has a coherent splitting
given by $\phi=e^{35}$ with $H^{0,2}=H^{0,3}=0$.
\item $A_{5,11}=\left(e^{15}+e^{25},e^{35}+e^{25},e^{35},- pe^{45} ,0\right)$  has a coherent splitting  given by $\phi=e^{35}$ with $H^{0,2}=H^{0,3}=0$.
\item $A_{5,12}=\left(e^{15}+e^{25},e^{25}+e^{35},e^{35}+e^{45},e^{45},0\right)$ has a coherent splitting  given by $\phi=e^{45}$ with $H^{0,2}=H^{0,3}=0$.
\item $A_{5,13}=\left(e^{15}, pe^{25} , q e^{35}+ re^{45} , q e^{45}- re^{35} ,0\right)$, where we assume $r\neq0$ (as $A_{5,13}(p,q,0)$ is isomorphic to $A_{5,7}(p,q,q)$), has a  coherent splitting given by $\phi=e^{25}$ with $H^{0,2}=H^{0,3}=0$ if $q\neq0,-\frac12$. If $q=-1/2$, $p\neq1$ then the same holds of $\phi=e^{15}$. The only cases left out are $(p,q)=(1,-\frac12)$ and $q=0$. In general, a basis of $Z^3$ is given by
\[ e^{145}, e^{245}, e^{135}, e^{345}, e^{125},e^{235} \]
 plus $e^{134}$ if $1+2q$ is zero, plus $e^{234}$ if $p+2q$ is zero. A basis of $Z^2$ is given by
\[e^{45} , e^{25}, e^{15}, e^{35}\]
plus $e^{34}$ if $q$ is zero, plus $e^{12}$ if $p=-1$. Thus, if $(p,q)=(1,-\frac12)$, then $(Z^2)^2=0$ contradicting \eqref{eqn:Z2wedgealpha} for any $\alpha$. On the other hand, if $q=0$, then $\alpha$ cannot be independent of $\beta=e^5$, as $\dim L_\beta(Z^3)<2$. But then  $\eqref{eqn:Z2wedgealpha}$ is only satisfied if $(p,q)=(-1,0)$, in which case we already know that a hypo structure exists.
 \item $A_{5,14}=\left(e^{25},0,e^{45}+ pe^{35} ,-e^{35}+ pe^{45} ,0\right)$ has a coherent splitting  given by $\phi=e^{25}$ with $H^{0,2}=H^{0,3}=0$ if $p\neq0$.
\item $A_{5,15}=\left(e^{15}+e^{25},e^{25},e^{45}+ pe^{35} , pe^{45} ,0\right)$ has a coherent splitting  given by $\phi=e^{45}$ with $H^{0,2}=H^{0,3}=0$ if $p\neq-1$.
\item $A_{5,16}=\left(e^{25}+e^{15},e^{25}, pe^{35} + qe^{45} , pe^{45} - qe^{35} ,0\right)$ has a coherent splitting  given by $\phi=e^{25}$ if $p\neq0$. If $p=0$, then $\beta=e^5$ and $L_\beta(Z^3)=0$.
\item $A_{5,17}=\left( p e^{15}+e^{25}, p e^{25}-e^{15}, r e^{45}+ q e^{35},- r e^{35}+ q e^{45},0\right)$, $r\neq 0$. Then $Z^3\wedge e^5=0$. So
if $p+q\neq0$, then $\beta=e^5$ and Proposition~\ref{prop:unimodular} applies. Otherwise, the same argument together with Proposition~\ref{prop:third} shows that necessarily $\alpha=e^5$. Now if $p+q=0$ but $p\neq 0$ and $r\neq\pm1$, then $Z^2=e^5\wedge\Lambda^1$, so by \eqref{eqn:Z2wedgealpha} no hypo structure exists.
\item $A_{5,18}=\left( p e^{15}+e^{35}+e^{25}, p e^{25}+e^{45}-e^{15}, p e^{35}+e^{45},- p e^{45}-e^{35},0\right)$. Then $\beta=pe^5$, and $L_{e^5}(Z^3)=0$. So if $p\neq 0$ we obtain an obstruction.
\end{itemize}
\end{proof}
\begin{table}
\caption{\label{table:IrreducibleNonContact}Nondecomposable, non-contact hypo Lie algebras}
\[
 \begin{array}{|l|l|l|}
 \hline
 \text{Name} & \text{Structure constants} & \text{Hypo structure}\\
\hline
A_{5,1}&\left(e^{35},e^{45},0,0,0\right) & e^1,e^3,e^2,e^4,e^5\\
 A_{5,2}&\left(e^{25},e^{35},e^{45},0,0\right) & e^1,e^4,e^3,e^2,e^5\\
 A_{5,7}(p,-p,-1)&\left(e^{15}, pe^{25} ,- pe^{35},- e^{45},0\right)&e^1,e^4,e^2,e^3,e^5\\
A_{5,8}(-1)&\left(e^{25},0,e^{35}, -e^{45} ,0\right)&e^1,e^2,e^3,e^4,e^5\\
A_{5,13}(-1,0,r)&\left(e^{15},- e^{25}, re^{45} ,- re^{35} ,0]\right)&e^{1},e^{2},e^{3},e^{4},e^{5}\\
A_{5,14}(0)&\left(e^{25},0,e^{45},- e^{35},0\right)&  e^{1}, e^{2},  e^{3},  e^{4},  e^{5}\\
A_{5,15}(-1)&\left(e^{15}+e^{25},e^{25},e^{45}-e^{35},- e^{45},0\right)&e^{1},e^{4},e^{3},e^{2},e^{5}\\
A_{5,17}(0,0,r)&\left(e^{25},- e^{15}, re^{45} ,- e^{35} r,0\right)&e^{1},e^{2},e^{3},e^{4},e^{5}\\
A_{5,17}(p,-p,1)&\left(e^{25}+ p e^{15},-e^{15}+ p e^{25},e^{45}- p e^{35},-e^{35}- p e^{45},0\right)&e^{1},e^{3},e^{2},e^{4},e^{5}\\
A_{5,17}(p,-p,-1)&\left(e^{25}+ p e^{15},-e^{15}+ p e^{25},-e^{45}- p e^{35},e^{35}- p e^{45},0\right)&e^{1},e^{3},e^{4},e^{2},e^{5}\\
A_{5,18}(0)&\left(e^{35}+e^{25},-e^{15}+e^{45},e^{45},- e^{35},0\right)&e^{1},e^{3},e^{2},e^{4},e^{5}\\
\hline
 \end{array}\]
\end{table}

\section{Decomposable contact Lie algebras}
By \cite{Diatta}, there are two types of decomposable $5$-dimensional  Lie algebras with an invariant contact form.
First, the Lie algebras $(0,e^{12})\oplus \lie{g}_3$, where $(0,e^{12})$ is the Lie algebra of affine transformations of $\R$, and $\lie{g}_3$ is any Lie algebra of dimension three other than
$(0,e^{12},e^{13})$ or $(0,0,0)$. Second, the Lie algebras of the form $\lie{g}_4\oplus\R$, where $\lie{g}_4$ is a four-dimensional Lie algebra carrying an exact symplectic form. In this section we show that only one of these families admits a hypo structure, and it belongs to the first type.
\begin{proposition}
\label{prop:decomposable012}
If $\lie{g}_3$ is a solvable Lie algebra of dimension three,
then $(0,e^{12})\oplus \lie{g}_3$ has a hypo structure if and only if
$\lie{g}_3=A_{3,8}=(e^{23},- e^{13},0)$.
\end{proposition}
\begin{proof}
First, observe that
$A_{3,8}\oplus(0,e^{12})=(e^{23},- e^{13},0,0,e^{45})$
has a hypo structure given by the coframe
$e^{1},e^{2},e^{4},e^{3},e^{5}$.

To prove uniqueness, we resort once again to the list in \cite{Basarab-Horwath}. There are nine families of solvable Lie algebras of dimension three, of which the following five have a coherent splitting with $H^{0,2}=H^{0,3}=0$:
\begin{align*}
A_{3,2}\oplus(0,e^{12})&=(0,e^{12},0,0,e^{45}),& \phi&=e^{14}\\
A_{3,4}\oplus(0,e^{12})&=(e^{23}+e^{13},e^{23},0,0,e^{45}), & \phi&=e^{34}\\
A_{3,5}\oplus(0,e^{12})&=(e^{13},e^{23},0,0,e^{45}),& \phi&=e^{34}\\
A_{3,7}\oplus(0,e^{12})&=(e^{13}, q e^{23},0,0,e^{45}), \quad 0<\abs{q}<1, & \phi&=e^{34}\\
A_{3,9}\oplus(0,e^{12})&=(q e^{13}+e^{23}, q e^{23}-e^{13},0,0,e^{45}), \quad q>0, & \phi&=e^{34}.
\end{align*}
Therefore, by Proposition \ref{prop:criterion}, there is no hypo structure on
$A_{3,i}$, for $i=2, 4, 5, 7, 9$.

For $A_{3,1}\oplus(0,e^{12})=(0,0,0,0,e^{45})$, we find that $\beta=e^{4}$ and
$L_\beta(Z^3)$ is spanned by $e^{1234}$, so Proposition~\ref{prop:unimodular} applies.

For $A_{3,3}\oplus(0,e^{12})=(e^{23},0,0,0,e^{45})$, we have that $\beta=- e^{4}$ and
\[L_\beta(Z^3)  =\operatorname{Span}\{- e^{1234},e^{2345}\}.\]
So, by Proposition~\ref{prop:unimodular}, if $(\alpha,\omega_i)$ is a hypo structure then
\[\alpha\in \operatorname{Span}\{e^{2},e^{3},e^{4}\},\]
implying that $Z^2\wedge\alpha\wedge\beta$ is contained in $L_\beta(Z^3)$, which is absurd.

Finally, the Lie algebra
$A_{3,6}\oplus(0,e^{12})=(e^{13},- e^{23},0,0,e^{45})$ satisfies
\begin{gather*}
Z^2=\Span{e^{12},e^{13},e^{23},e^{34},e^{45}},\\
Z^3=\Span{e^{123},e^{124},e^{134},e^{145}-e^{135},e^{234},e^{245}+e^{235},e^{345}}.
\end{gather*}
So $\alpha$ lies in the span of $e^3,e^4$ by Lemma~\ref{lemma:alphaXzero}. Moreover $\beta=e^4$, thus $\alpha$ has the form $e^3+ae^4$. Defining $V$ as in Proposition~\ref{prop:third}, we see that $e^4\wedge V$ is contained in the span of $e^{1234}$ and $e^{2345}$. Since $Z_2\wedge e^{34} = e^{1234}$, there is no hypo structure.
\end{proof}
\begin{remark}
Notice that Proposition~\ref{prop:decomposable012} does not  apply to contact Lie algebras alone, but also to the non-contact Lie algebras
$A_{3,1}\oplus(0,e^{12})$ and $A_{3,5}\oplus(0,e^{12})$.
\end{remark}
Decomposable contact Lie algebras of the type $\lie{g}_4\oplus\R$ are
not unimodular, because the volume form is exact, and so it makes sense to apply Proposition~\ref{prop:unimodular}. This turns out to be sufficient in order to show that no hypo structure exists on these Lie algebras.
\begin{proposition}
\label{prop:exactsymplectic}
If $\lie{g}_4$ is a $4$-dimensional solvable Lie algebra with an exact
symplectic form, then $\lie{g}_4\oplus\R$ has no hypo structure.
\end{proposition}
\begin{proof}
Observe that $\lie{g}_4$ is necessarily indecomposable,
because it admits an exact symplectic form. From the list in \cite{Basarab-Horwath}, $\lie{g}_4$ must belong to one of four families, to each of which we apply Proposition~\ref{prop:unimodular}:
\begin{itemize}
\item The Lie algebra $A_{4,7}\oplus\R=(e^{23}+2e^{14},e^{24}+e^{34},e^{34},0,0)$ has $\beta=e^4$, and $L_\beta(Z^3)$ has dimension one.
\item The Lie algebra $A_{4,8}\oplus\R=(e^{23}+(1+q)e^{14},e^{24},qe^{34},0,0)$, $-1< q\leq 1$ has $\beta=e^4$, and $L_\beta(Z^3)$ has dimension one except if $q=-1/2$. In this case, it is spanned by $e^{1345}, e^{2345}$, and
\[Z^2=\Span{e^{13},e^{23}+\frac{1}{2} e^{14},e^{24},e^{34},e^{45}},\]
so by Proposition~\ref{prop:unimodular} there is no hypo structure.
\item The Lie algebra
$A_{4,9}\oplus\R=\left(e^{23}+2qe^{14},qe^{24}+e^{34},-e^{24}+qe^{34},0,0\right)$,
with $q> 0$. Then $\beta=e^4$, and $L_\beta(Z^3)$ is one-dimensional.
\item The Lie algebra $A_{4,10}\oplus\R=(e^{13}+e^{24},e^{23}-e^{14},0,0,0)$ has $\beta=e^3$,
and \[L_\beta(Z^3)=\Span{e^{2345},e^{1345}}.\] So $L_\alpha$ kills $L_\beta(Z^3)$ if and only if $\alpha$ lies in the span of $e^3,e^4,e^5$, in which case $\alpha\wedge\beta\wedge Z^2$ is contained in $L_\beta(Z^3)$.
\end{itemize}
Thus, in neither case is there a hypo structure.
\end{proof}
\begin{remark}
In the proof of Proposition~\ref{prop:exactsymplectic}, we have left out the cases $A_{4,8}(-1)$ and $A_{4,9}(0)$ because they do not have any exact symplectic form. They  have no hypo structure either. Indeed,  $A_{4,8}(-1)$ has a coherent splitting $\phi=e^{24}$ with $H^{0,2}=H^{0,3}=0$. For  $A_{4,9}(0)$, we apply Lemma~\ref{lemma:alphaXzero} to show that $\alpha$ has no component along $e^1$, contradicting \eqref{eqn:Z2wedgealpha}.
\end{remark}

\section{Decomposable non-contact Lie algebras}
Decomposable Lie algebras of dimension five may either be of the form $\lie{g}_3\oplus \lie{h}_2$, where we are allowing the factors themselves to be decomposable, or $\lie{g}_4\oplus\R$.  In the former case, by Proposition~\ref{prop:decomposable012} we can assume $\lie{h}_2=\R^2$.  Without resorting to Mubarakzyanov's classification, we can characterize which of these Lie algebras have a hypo structure.
\begin{proposition}
 Let $\lie{g}_3$ be a Lie algebra of dimension $3$. Then $\lie{g}=\lie{g}_3\oplus\R^2$ admits
a hypo structure if and only if $\lie{g}_3$ is unimodular.
\end{proposition}
\begin{proof}
Let $e^1,\dotsc, e^5$ be a coframe reflecting the splitting $\lie{g}=\lie{g}_3\oplus\R^2$, so that $e^1,e^2,e^3$ is a basis of $\lie{g}_3^*\subset\lie{g}^*$ and $e^4,e^5$ a basis of $(\R^2)^*$.

If $\lie{g}_3$ is unimodular, then the coframe $e^1,e^2,e^4,e^5,e^3$ determines a hypo structure by \eqref{eqn:adaptedframe}, because $e^{12},e^{13},e^{23}$ are closed.

If $\lie{g}_3$ is not unimodular, we can assume that $e^3=\beta$ as defined in Lemma~\ref{lemma:unimodular}. Then $e^3$ is closed and $de^{12}\neq 0$. Moreover $de^i\wedge e^3=0$,
$i=1,2$. This is because if $e_1,e_2,e_3$ is a basis of $\lie{g}_3$ dual to $e^1,e^2,e^3$, then
\[0=\operatorname{tr}(\ad e_2)=e^1([e_2,e_1])+e^3([e_2,e_3])=e^1([e_2,e_1]).\]
Thus $Z^3=(e^3\wedge\Lambda^2)\oplus W$, where 	\[W\subset\Lambda^3(\Span{e^1,e^2,e^4,e^5}).\] Since $de^{12}\neq 0$, $W\subset \Span{e^{145},e^{245}}$;
so $L_\beta(Z_3)$ has the same dimension as $\ker d\cap \Span{e^1,e^2}$, which is at most one since $\lie{g}_3$ is not abelian. By Proposition~\ref{prop:unimodular} this concludes the proof.
\end{proof}

For the other case, we must refer to Mubarakzyanov's classification.
\begin{proposition}
If $\lie{g}_4$ is an indecomposable solvable Lie algebra of dimension four, then $\lie{g}_4\oplus\R$ has no hypo structure.
\end{proposition}
\begin{proof}
By  \cite{Basarab-Horwath}, there are $10$ families $A_{4,1},\dotsc, A_{4,10}$ of solvable Lie algebras of dimension four.
The families $A_{4,7}$ through $A_{4,10}$ have no hypo structure by Proposition~\ref{prop:exactsymplectic} and the subsequent remark. The following families have a coherent splitting with $H^{0,2}=H^{0,3}=0$:
\begin{align*}
A_{4,1}\oplus\R&=(e^{24},e^{34},0,0,0), & \phi&=e^{34}\\
A_{4,2}\oplus\R&=(q e^{14} ,e^{24}+e^{34},e^{34},0,0), \quad q\neq0, & \phi&=e^{34}\\
A_{4,3}\oplus\R&=\left(e^{14},e^{34},0,0,0\right),& \phi&=e^{34}\\
A_{4,4}\oplus\R&=(e^{14}+e^{24},e^{24}+e^{34},e^{34},0,0), & \phi&=e^{45}\\
A_{4,5}\oplus\R&=\left(e^{14}, q e^{24}, p e^{34},0,0\right), \quad \quad p,q\neq0,
& \phi&=\begin{cases} e^{14}, &(p,q)=(-1,-1) \\ e^{34}, &q\neq-1 \\ e^{24}, &p\neq -1\end{cases}
\end{align*}
Notice that  $A_{4,1}\oplus\R$ is nilpotent and isomorphic to $(0,0,0,e^{12},e^{14})$.

The remaining family is  \[A_{4,6}\oplus\R=\left(q e^{14},e^{34}+ p e^{24} ,-e^{24}+ p e^{34},0,0\right), \quad q\neq 0, \quad p\geq 0.\]
If $p>0$, $\phi=e^{14}$ defines a coherent splitting with $H^{0,2}=H^{0,3}=0$. If $p=0$, then
$\beta=e^4$, and
\[Z^3=\Span{ e^{235}, e^{145} ,  e^{245},e^{134}, e^{234}, e^{124}, e^{345}}.\]
Therefore $L_\beta(Z^3)$ has dimension one, and Proposition~\ref{prop:unimodular} applies.
\end{proof}

Applying the classification of three-dimensional Lie algebras, we finally obtain:
\begin{theorem}
A solvable Lie algebra of dimension five has a hypo structure if and only if it appears in the list of Proposition~\ref{prop:diattahypo},
  it appears in Table~\ref{table:IrreducibleNonContact}, or
 it is one of the following:
\begin{multline*}
 (0,0,0,0,0), \quad  (e^{23},0,0,0,0), \quad (e^{23},-e^{13},0,0,0), \quad  (e^{13},-e^{23},0,0,0), \\
 (e^{23},-e^{13},0,0,e^{45}).
\end{multline*}
\end{theorem}

\smallskip

{\small

\small\noindent Dipartimento di Matematica e Applicazioni, Universit\`a di Milano Bicocca,  Via Cozzi 53, 20125 Milano, Italy.\\
\texttt{diego.conti@unimib.it}

\smallskip
\small\noindent Universidad del Pa\'{\i}s Vasco, Facultad de Ciencia y Tecnolog\'{\i}a, Departamento de Matem\'aticas,
Apartado 644, 48080 Bilbao, Spain. \\
\texttt{marisa.fernandez@ehu.es}\\
\texttt{joseba.santisteban@ehu.es}
\end{document}